\documentclass[a4paper,12pt,intlimits,oneside]{amsart}

\usepackage{enumerate}

\usepackage{ latexsym,amssymb}

\newcommand{\comment}[1]{}

\newcommand{\epf}{ $\Box$\medskip}

\numberwithin{equation}{section}

\def\lsim{\raisebox{-1ex}{$~\stackrel{\textstyle <}{\sim}~$}}

\newcounter{rea}
\newcounter{rej}
\newcounter{res}
\newtheorem{thm}{Theorem}[section]

\newtheorem{cor}[thm]{Corollary}
\newtheorem{lem}[thm]{Lemma}

\begin{document}

\title[]{Intrinsic square functions on functions spaces including weighted Morrey spaces}
\author[J. Feuto]{Justin Feuto}
\address{Laboratoire de Math\'ematiques Fondamentales, UFR Math\'ematiques et Informatique, Universit\'e de Cocody, 22 B.P 1194 Abidjan 22. C\^ote d'Ivoire}
\email{{\tt justfeuto@yahoo.fr}}

\subjclass{42B25; 47B38;  47B47}
\keywords{  amalgams spaces, Morrey spaces, commutator, $g$-function of Littlewood-Paley, Lusin area function.}
\date{}




\begin{abstract}
We prove that the intrinsic square functions including Lusin area integral and Littlewood-Paley $g^{\ast}_{\lambda}$-function as defined by Wilson, are bounded in a class of function spaces include weighted Morrey spaces. The corresponding commutators generated by $BMO$ functions are also considered. 
\end{abstract}


\maketitle 

\section{Introduction}
The classical Morrey spaces were introduced by Morrey \cite{Mo} in connection with partial differential equations. We recall that a real-valued function $f$ is said to belong to the Morrey space $L^{q,\lambda}$ on the $n$-dimensional euclidean space $\mathbb R^{n}$ provided the following norm is finite:
\begin{equation*}
\left\|f\right\|_{L^{q,\lambda}}:=\left(\sup_{(y,r)\in\mathbb R^{n}\times\mathbb R^{\ast}_{+}}r^{\lambda-n}\int_{B(y,r)}\left|f(x)\right|^{q}dx\right)^{\frac{1}{q}}.
\end{equation*}
Here $1\leq q<\infty$, $0<\lambda<n$, $\mathbb R^{\ast}_{+}=\left(0,\infty\right)$ and $B(y,r)$ is a ball in $\mathbb R^{n}$ centered at $y$ of radius $r$.

Chiarenza and Frasca \cite{CF} established the boundedness of the Hardy-Littlewood maximal operator, the fractional operator and Calder\'on-Zygmund operator on these spaces. These operator are also bounded on Lebesgue spaces, and in weighted Lebesgue space  \cite{CFe,Mu}.

Twenty years ago, Fofana introduced a class of function spaces comprising Lebesgue and Morrey spaces \cite{Fo0}. Precisely, for $1\leq q\leq p\leq\infty$, let $(L^{q},L^{p})(\mathbb R^{n})$ be the Wiener amalgam space of $L^{q}(\mathbb R^{n})$ and $L^{p}(\mathbb R^{n})$, i.e., the space  of measurable functions $f:\mathbb R^{n}\rightarrow\mathbb C$ which are locally in $L^{q}(\mathbb R^{n})$ and such that the function  $y\mapsto\left\|f\chi_{B(y,1)}\right\|_{q}$ belongs to $L^{p}(\mathbb R^{n})$, where for $r>0$, $B(y,r)=\left\{x\in\mathbb R^{n}/\left|x-y\right|<r\right\}$ is the open ball centered at $y$ with radius $r$, $\chi_{B(y,r)}$ its characteristic function and $\left\|\cdot\right\|_{q}$ denoting the usual Lebesgue norm in $L^{q}(\mathbb R^{n})$. As we can see in \cite{H,FS}, we have the following properties. 
\begin{itemize}
\item For $1\leq q\leq p\leq\infty$, the space $(L^{q},L^{p})(\mathbb R^{n})$ is a Banach space when it is equipped with the norm 
\begin{equation*}
\left\|f\right\|_{q,p}:=\left(\int_{\mathbb R^{n}}\left\|f\chi_{B(y,1)}\right\|^{p}_{q}\right)^{\frac{1}{p}}
\end{equation*}
with the usual modification when $p=\infty$.
\item The amalgam space $(L^{q},L^{q})$ is equal to the Lebesgue space $L^{q}$ with equivalence norms, provided $q=p$, while for  $q\leq s\leq p$, we have $L^{s}(\mathbb R^{n})$ continuously embedded in $(L^{q},L^{p})(\mathbb R^{n})$.
\end{itemize}

In Lebesgue spaces $L^{q}(\mathbb R^{n})$, it is well known that for $r>0$ and $x\in\mathbb R^{n}$, the dilation operators $\delta^{q}_{r}:f\mapsto r^{\frac{n}{q}}f(r\cdot)$ and the translation operators $\tau_{x}:f\mapsto f(\cdot-x)$ are isometries. We use the usual convention that $\frac{1}{\infty}=0$. When we consider the amalgam spaces $(L^{q},L^{p})(\mathbb R^{n})$ with $q<p$, only translation operators conserve this property. 
 But it is easy to see that $f\in (L^{q},L^{p})$ if and only if we have 
\begin{equation*}
\left\|\delta^{\alpha}_{r}f\right\|_{q,p}<\infty,
\end{equation*}
for all $r>0$ and all $\alpha>0$. Notice that for $1\leq q,p,\alpha\leq\infty$, $r>0$ and $\alpha>0$, we have
\begin{equation}
\begin{array}{lll}
\left\|\delta^{\alpha}_{r}f\right\|_{q,p}&=&r^{n(\frac{1}{\alpha}-\frac{1}{q}-\frac{1}{p})}\left(\int_{\mathbb R^{n}}\left\|f\chi_{B(y,r)}\right\|^{p}_{q}dy\right)^{\frac{1}{p}}\\
&\approx& \left[\int_{\mathbb R^{n}}\left(\left|B(y,r)\right|^{\frac{1}{\alpha}-\frac{1}{q}-\frac{1}{p}}\left\|f\chi_{B(y,r)}\right\|_{q}\right)^{p}dy\right]^{\frac{1}{p}},
\end{array}\footnote{Hereafter we propose the following abbreviation $\mathrm{\bf A}\approx \mathrm{\bf B}$ for the inequalities $C^{-1}\mathrm{\bf A}\leq\mathrm{\bf B}\leq C\mathrm{\bf A}$, where $C$ is a positive constant independent of the main parameters.}\label{equivamalg}
\end{equation}
where $\left|B(y,r)\right|$ stands for the Lebesgue measure of the ball $B(y,r)$.
This bring Fofana \cite{Fo0} to consider the subspace $(L^{q},L^{p})^{\alpha}(\mathbb R^{n})$ of $(L^{q},L^{p})(\mathbb R^{n})$ that consists in measurable functions $f$ such that $\left\|f\right\|_{q,p,\alpha}<\infty$, where for $1\leq q,p,\alpha\leq\infty$,
\begin{equation*}
\left\|f\right\|_{q,p,\alpha}:=\sup_{r>0}\left\|\delta^{\alpha}_{r}f\right\|_{q,p}.
\end{equation*}

As proved in \cite{Fo0,FFK1,FFK2}, the spaces $(L^{q},L^{p})^{\alpha}(\mathbb R^{n})$ are non trivial if and only if $q\leq\alpha\leq p$. In this case, for fixed $1\leq q<\alpha$  and $p$ varying from $\alpha$ to $\infty$, these spaces form a chain of distinct Banach spaces beginning with Lebesgue space $L^{\alpha}(\mathbb R^{n})$ and ending by the classical Morrey space $L^{q,\frac{nq}{\alpha}}(\mathbb R^{n})=(L^{q},L^{\infty})^{\alpha}(\mathbb R^{n})$. More precisely, we have the following continuous injections 
 $$L^{\alpha}(\mathbb R^{n})\hookrightarrow (L^{q},L^{p_{1}})^{\alpha}(\mathbb R^{n})\hookrightarrow (L^{q},L^{p_{2}})^{\alpha}(\mathbb R^{n})\hookrightarrow L^{q,\frac{nq}{\alpha}}(\mathbb R^{n}),$$
 for $q\leq \alpha<p_{1}<p_{2}<\infty$. 
It is therefore interesting to know the behavior of operators which are bounded on Lebesgue and Morrey spaces, on these spaces. 
%
%

We proved in \cite{BFF} that classical operators such as the Hardy-Littlewood maximal operator, the Calder\'on-Zygmund operator and Riesz potentials which are known to be bounded in Lebesgue and Morrey spaces, are also bounded $(L^{q},L^{p})^{^{\alpha}}$ spaces (see also \cite{FFK3}). 

In \cite{KS}, Komori and Shirai considered the weighted Morrey spaces $L^{q,\kappa}_{w}(\mathbb R^{n})$ when studying the boundedness of Hardy-Littlewood and Calder\'on-Zygmund operators. 

Let $0<\kappa<1$ and $w$ a weight on $\mathbb R^{n}$, i.e., a positive locally integrable function on $\mathbb R^{n}$. The weighted Morrey space $L^{q,\kappa}_{w}(\mathbb R^{n})$, consists of measurable functions $f$ such that $\left\|f\right\|_{L^{q,\kappa}_{w}}<\infty$, where
\begin{equation*}
\left\|f\right\|_{L^{q,\kappa}_{w}}:=\sup_{B}\left(\frac{1}{w(B)^{\kappa}}\int_{B}\left|f(x)\right|^{q}w(x)dx\right)^{\frac{1}{q}}.
\end{equation*}
These spaces generalize weighted Lebesgue spaces $L^{q}_{w}(\mathbb R^{n})$, namely the space consisting in measurable functions $f$ satisfying 
\begin{equation*}
\left\|f\right\|_{q_{w}}:=\left(\int_{\mathbb R^{n}}\left|f(x)\right|^{q}w(x)dx\right)^{\frac{1}{q}}<\infty.
\end{equation*}
  In this work, we consider for $1\leq q\leq\alpha\leq p\leq\infty$ and a weight $w$, the space $(L^{q}_{w},L^{p})^{\alpha}(\mathbb R^{n})$ consists of measurable functions $f$ such that $\left\|f\right\|_{q_{w},p,\alpha}<\infty$, where 
\begin{equation*}
\ _{r}\left\|f\right\|_{q_{w},p,\alpha}:=\left[\int_{\mathbb R^{n}}\left(w(B(y,r))^{\frac{1}{\alpha}-\frac{1}{q}-\frac{1}{p}}\left\|f\chi_{B(y,r)}\right\|_{q_{w}}\right)^{p}dy\right]^{\frac{1}{p}},
\end{equation*}
for $r>0$, and 
\begin{equation*}
\left\|f\right\|_{q_{w},p,\alpha}:=\sup_{r>0}\ _{r}\left\|f\right\|_{q_{w},p,\alpha}, 
\end{equation*}
with $w(B(y,r))=\int_{B(y,r)}w(x)dx$ and the usual modification when $p=\infty$.
 When $w\equiv 1$, we recover $(L^{q},L^{p})^{\alpha}(\mathbb R^{n})$ spaces while for $q<\alpha$ and $p=\infty$, the spaces $(L^{q}_{w},L^{\infty})^{\alpha}(\mathbb R^{n})$ are noting but the weighted Morrey spaces $L^{q,\kappa}_{w}(\mathbb R^{n})$, with $\kappa=\frac{1}{q}-\frac{1}{\alpha}$. Wilson in \cite{Wi2} proved that for $1<q<\infty$ and $0<\gamma\leq 1$, the intrinsic square operators $S_{\gamma}$ given by Relation (\ref{intsqf}), are bounded in the weighted Lebesgue spaces $L^{q}_{w}$, 
whenever the weight $w$ fulfilled the $\mathcal A_{q}$ condition of Muckenhoupt. Wang extend this result to weighted Morrey spaces $L^{q,\kappa}_{w}(\mathbb R^{n})$. We prove here that these operators and others known to be bounded on weighted Lebesgue and weighted Morrey spaces, are also bounded in the more general setting of $(L^{q}_{w},L^{p})^{\alpha}$ spaces. 
%
 
 This paper is organized as follows:
 
 In the second section, we recall the definitions of the operators we are going to consider and recall the results on weighted Lebesgue and Morrey spaces. In  Section 3 we state our results and in the last section we give their proofs.
 
 Throughout the paper, the letter $C$ is used for non-negative constants that may change from one occurrence to another. 
  The notation $\mathrm{\bf A}\lsim \mathrm{\bf B}$ will always stand for $\mathrm{\bf A}\leq C\mathrm{\bf B}$, where $C$ is a positive constant independent of the main parameters. 
 For $\alpha>0$ and a ball $B\subset\mathbb R^{n}$, we write $\alpha B$ for the ball with same center as $B$ and with radius $\alpha$ times radius of $B$. For any subset $E$ of $\mathbb R^{n}$, we denote $E^{c}:=\mathbb R^{n}\setminus E$ the complement of $E$. We denote by $\mathbb N^{\ast}$ the set of all positive integers.

\section{Definitions and known results}


For $0<\gamma\leq 1$, we denote by $\mathcal C_{\gamma}$ the family of function $\varphi$ defined on $\mathbb R^{n}$ with support in the closed unit ball $\mathbb B=\left\{x\in\mathbb R^{n}:\left|x\right|\leq 1\right\}$ and vanishing integral, i.e., $\int_{\mathbb R^{n}}\varphi(x)dx=0$, and such that for all $x,x'\in\mathbb R^{n}$, $\left|\varphi(x)-\varphi(x')\right|\leq\left|x-x'\right|^{\gamma}$. Let $\mathbb R^{n+1}_{+}=\mathbb R^{n}\times\left(0,\infty\right)$ and $\varphi_{t}(x)=t^{-n}\varphi(t^{-1}x)$. 
The intrinsic square function of $f$ (of order $\gamma$) is defined by the formula
\begin{equation}
S_{\gamma}(f)(x)=\left[\int_{\Gamma(x)}\left(\sup_{\varphi\in\mathcal C_{\gamma}}\left|f\ast\varphi_{t}(y)\right|\right)^{2}\frac{dydt}{t^{n+1}}\right]^{\frac{1}{2}},\label{intsqf}
\end{equation}
where for $x\in\mathbb R^{n}$, $\Gamma(x)$ denote the usual "cone of arperture one", $$\Gamma(x)=\left\{(y,t)\in\mathbb R^{n+1}_{+}:\left|x-y\right|<t\right\}.$$
For $1<q<\infty$ and $0<\gamma\leq 1$, the operators $S_{\gamma}$ are bounded on $L^{q}_{w}(\mathbb R^{n})$ provided $w\in\mathcal A_{q}$ \cite{Wi}. 
We recall that a weight $w$ is of class $\mathcal A_{q}$ or belongs to $\mathcal A_{q}$ for $1<q<\infty$ if there exists a constant $C>0$ such that for all balls $B\subset\mathbb R^{n}$ we have
\begin{equation}
\left(\frac{1}{\left|B\right|}\int_{B}w(x)dx\right)\left(\frac{1}{\left|B\right|}\int_{B}w^{\frac{-1}{q-1}}(x)dx\right)^{q-1}\leq C.\label{apcondition}
\end{equation}
In the setting of weighted Morrey spaces one has the following.
\begin{thm}[Theorem 1.1 \cite{Wa}]\label{rap1}
Let $0<\gamma\leq 1$, $1<q<\infty$, $0<\kappa<1$ and $w\in\mathcal A_{q}$. Then there exists $C>0$ such that 
\begin{equation}
\left\|S_{\gamma}f\right\|_{L^{q,\kappa}_{w}}\leq C\left\|f\right\|_{L^{q,\kappa}_{w}}.
\end{equation}
\end{thm}

 We also define the intrinsic Littlewood-Paley $g$-function $g_{\gamma}(f)$ and $g^{\ast}_{\lambda}$-function $g^{\ast}_{\lambda,\gamma}(f)$ by 
\begin{equation*}
g_{\gamma}(f)(x)=\left[\int^{\infty}_{0}\left(\sup_{\varphi\in\mathcal C_{\gamma}}\left|f\ast\varphi_{t}(y)\right|\right)^{2}\frac{dt}{t}\right]^{\frac{1}{2}},
\end{equation*}
and
\begin{equation*}
g^{\ast}_{\lambda,\gamma}(f)(x)=\left[\int_{\mathbb R^{n+1}_{+}}\left(\frac{t}{t+\left|x-y\right|}\right)^{\lambda n}\left(\sup_{\varphi\in\mathcal C_{\gamma}}\left|f\ast\varphi_{t}(y)\right|\right)^{2}\frac{dydt}{t^{n+1}}\right]^{\frac{1}{2}},
\end{equation*}
respectively.
\begin{thm}[Theorem 1.3 \cite{Wa}]\label{rap2}
Let $0<\gamma\leq 1$, $1<q<\infty$, $0<\kappa<1$ and $w\in\mathcal A_{q}$. If $\lambda>\max\left\{q,3\right\}$, then there exists $C>0$ such that 
\begin{equation}
\left\|g^{\ast}_{\lambda,\gamma}f\right\|_{L^{q,\kappa}_{w}}\leq C\left\|f\right\|_{L^{q,\kappa}_{w}}.
\end{equation}
\end{thm}
Let $b$ be a locally integrable function. The commutator of $b$ and $S_{\gamma}$ is defined by
\begin{equation*}
\left[b,S_{\gamma}\right](f)(x)=\left(\int_{\Gamma(x)}\sup_{\varphi\in\mathcal C_{\gamma}}\left|\int_{\mathbb R^{n}}(b(x)-b(z))\varphi_{t}(y-z)f(z)dz\right|^{2}\frac{dydt}{t^{n+1}}\right)^{\frac{1}{2}},
\end{equation*} 
and the commutator of $b$ and $g^{\ast}_{\lambda,\gamma}$ by
\begin{equation*}
\left[b,g^{\ast}_{\lambda,\gamma}\right](f)(x)=\left(\int_{\mathbb R^{n+1}_{+}}\left(\frac{t}{t+\left|x-y\right|}\right)^{\lambda n}\sup_{\varphi\in\mathcal C_{\gamma}}\left|\int_{\mathbb R^{n}}(b(x)-b(z))\varphi_{t}(y-z)f(z)dz\right|^{2}\frac{dydt}{t^{n+1}}\right)^{\frac{1}{2}}.
\end{equation*}
A locally integrable function $b$ belongs to $BMO(\mathbb R^{n})$ (bounded mean oscillation functions) if $\left\|b\right\|_{BMO(\mathbb R^{n})}<\infty$, where
$$\left\|b\right\|_{BMO(\mathbb R^{n})}:=\sup_{B:\text{ ball}}\frac{1}{\left|B\right|}\int_{B}\left|b(x)-b_{B}\right|dx.$$
We have the following result in the context of weighted Lebesgue spaces.
\begin{thm}[Theorem 3.1 \cite{Wa}]
Let $0<\gamma\leq 1$, $1<q<\infty$ and $w\in\mathcal A_{q}$. Then the commutators $\left[b,S_{\gamma}\right]$ and $\left[b,g^{\ast}_{\lambda,\gamma}\right]$ are bounded on $L^{q}_{w}(\mathbb R^{n})$ whenever $b\in BMO(\mathbb R^{n})$.
\end{thm}
For weighted Morrey spaces, the following results are proved.
\begin{thm}[Theorem 1.2 \cite{Wa}]\label{rap3}
Let $0<\gamma\leq 1$, $1<q<\infty$, $0<\kappa<1$ and $w\in\mathcal A_{q}$. Suppose that $b\in BMO$, then there exists $C>0$ such that
\begin{equation*}
\left\|\left[b,S_{\gamma}\right]f\right\|_{L^{q,\kappa}_{w}}\leq C\left\|f\right\|_{L^{q,\kappa}_{w}}.
\end{equation*}
\end{thm}
\begin{thm}[Theorem 1.4 \cite{Wa}]\label{rap4} Let $0<\gamma\leq 1$, $1<q<\infty$, $0<\kappa<1$ and $w\in\mathcal A_{q}$. If $b\in BMO(\mathbb R^{n})$ and $\lambda>\max\left\{q,3\right\}$, then there is a constant $C>0$ independent of $f$ such that 
\begin{equation*}
\left\|\left[b,g^{\ast}_{\lambda,\gamma}\right](f)\right\|_{L^{q,\kappa}_{w}}\leq C\left\|f\right\|_{L^{q,\kappa}_{w}}.
\end{equation*}
\end{thm}
\section{Statement of our main results}
Since our space at least for the case where the weight is
equal to 1, are included in Morrey spaces, we already know that the image is
in the space of Morrey. But what is shown is that if one has a slightly stronger
assumption, then this is also true for the image.

 For the intrinsic square function $S_{\gamma}$, we have the following result.
\begin{thm}\label{main1}
Let $0<\gamma\leq 1$, $1< q\leq\alpha< p\leq\infty$ and $w\in\mathcal A_{q}$. The operators $S_{\gamma}$ are bounded in $(L^{q}_{w},L^{p})^{\alpha}(\mathbb R^{n})$. 
\end{thm}
Theorem \ref{rap1} is a particular case of this result. The next, concerning the intrinsic Littlewood-Paley $g^{\ast}_{\lambda}$-function is an extension of Theorem \ref{rap2}.
\begin{thm}\label{main3}
Let $0<\gamma\leq 1$, $1<q\leq\alpha< p\leq\infty$ and $w\in\mathcal A_{q}$. If $\lambda>\max\left\{q,3\right\}$ then there exists a constant $C>0$ such that
\begin{equation*}
\left\|g^{\ast}_{\lambda,\gamma}(f)\right\|_{q_{w},p,\alpha}\leq C\left\|f\right\|_{q_{w},p,\alpha},
\end{equation*}
for all $f\in (L^{q}_{w},L^{p})^{\alpha}(\mathbb R^{n})$.
\end{thm}
As far as their commutators, we have the following results which are extensions of Theorems \ref{rap3} and \ref{rap4} respectively.
\begin{thm}\label{main2}
Let $0<\gamma\leq 1$, $1<q\leq\alpha< p\leq \infty$ and $w\in\mathcal A_{q}$. Suppose that $b\in BMO(\mathbb R^{n})$, then there exists a constant $C>0$ not depending on $f$ such that 
\begin{equation*}
\left\|\left[b,S_{\gamma}\right](f)\right\|_{q_{w},p,\alpha}\leq C\left\|f\right\|_{q_{w},p,\alpha},
\end{equation*}
for all $f\in (L^{q}_{w},L^{p})^{\alpha}(\mathbb R^{n})$.
\end{thm}
\begin{thm}\label{nouveau}
Let $0<\gamma\leq 1$, $1<q\leq\alpha<p\leq\infty$ and $w\in\mathcal A_{q}$. If $b\in BMO(\mathbb R^{n})$ and $\lambda>\max\left\{q,3\right\}$ then there exists a constant $C>0$ such that
\begin{equation*}
\left\|\left[b,g^{\ast}_{\lambda,\gamma}\right](f)\right\|_{q_{w},p,\alpha}\leq C\left\|f\right\|_{q_{w},p,\alpha},
\end{equation*}
for all $f\in (L^{q}_{w},L^{p})^{\alpha}(\mathbb R^{n})$.
\end{thm}

Since for any $0<\gamma\leq 1$ the functions $S_{\gamma}(f)$ and $g_{\gamma}(f)$ are pointwise comparable (see \cite{Wi}), as an immediate consequence of Theorems \ref{main1} and \ref{main2} we have the following results. 
\begin{cor}
Let $0<\gamma\leq 1$, $1< q\leq\alpha< p\leq\infty$ and $w\in\mathcal A_{q}$. The operator $g_{\gamma}$ is bounded in $(L^{q}_{w},L^{p})^{\alpha}(\mathbb R^{n})$.
\end{cor}
\begin{cor}
Let $0<\gamma\leq 1$, $1<q\leq\alpha< p\leq \infty$ and $w\in\mathcal A_{q}$. Suppose that $b\in BMO(\mathbb R^{n})$, then there exists a constant $C>0$ not depending on $f$ such that 
\begin{equation*}
\left\|\left[b,g_{\gamma}\right](f)\right\|_{q_{w},p,\alpha}\leq C\left\|f\right\|_{q_{w},p,\alpha},
\end{equation*}
for all $f\in (L^{q}_{w},L^{p})^{\alpha}(\mathbb R^{n})$.
\end{cor}
The above corollaries are extensions of Corollary 1.5 and Corollary 1.6 of \cite{Wa} respectively.
\section{Proof of the main results}
We will need the following properties of $\mathcal A_{q}$ weights (see Proposition 9.1.5 and Theorem 9.2.2 \cite{Gr}). Let $w\in \mathcal A_{q}$ for some $1<q<\infty$. 
\begin{enumerate}
\item For all $\lambda>1$ and all balls $B$, we have 
\begin{equation}
w(\lambda B)\lsim \lambda^{nq}w(B).\label{doubling}
\end{equation}
\item There exists a positive constant $\tau$ such that for every ball $B$, we have
\begin{equation}
\left(\frac{1}{\left|B\right|}\int_{B}w(t)^{1+\tau}dt\right)^{\frac{1}{1+\tau}}\lsim \frac{1}{\left|B\right|}\int_{B}w(t)dt,\label{reverseholder}
\end{equation}
and for any measurable subset $E$ of a ball $B$, we have
\begin{equation}
\frac{w(E)}{w(B)}\lsim \left(\frac{\left|E\right|}{\left|B\right|}\right)^{\frac{\tau}{1+\tau}}.\label{contwl}
\end{equation}
\end{enumerate} 
For our proofs, we use arguments as in \cite{FLY}.

\proof[Proof of Theorem \ref{main1}]

We fix $r>0$ and let $B=B(y,r)$ for some $y\in\mathbb R^{n}$. We write $f=f_{1}+f_{2}$, with $f_{1}=f\chi_{2B}$. Since $S_{\gamma}$ is a sublinear operator, we have 
\begin{equation} 
\left\|S_{\gamma}(f)\chi_{B}\right\|_{q_{w}}\leq\left\|S_{\gamma}(f_{1})\chi_{B}\right\|_{q_{w}}+\left\|S_{\gamma}(f_{2})\chi_{B}\right\|_{q_{w}}.\label{firstestim}
\end{equation}
For the term in $f_{1}$, we have 
\begin{equation}
\left\|S_{\gamma}(f_{1})\chi_{B}\right\|_{q_{w}}\lsim \left\|f\chi_{2B}\right\|_{q_{w}}\label{estf1}
\end{equation}
as an immediate consequence of the boundedness of $S_{\gamma}$ in $L^{q}_{w}(\mathbb R^{n})$. Our attention will be focused now on the second term. 

Let $\varphi\in\mathcal C_{\gamma}$, and $t>0$. Since the family $\mathcal C_{\gamma}$ is uniformly bounded with respect to the $L^{\infty}$-norm, we have 
\begin{equation}
\left|f_{2}\ast\varphi_{t}(u)\right|
\lsim t^{-n}\int_{(2B)^{c}\cap \tilde{B}(u,t)}\left|f(z)\right|dz,\label{estima1f2}
\end{equation}
for all $u\in\mathbb R^{n}$, where $\tilde{B}(u,t):=\left\{z\in\mathbb R^{n}/\left|z-u\right|\leq t\right\}$. 
Thus for all $x\in\mathbb R^{n}$, we have
\begin{eqnarray*}
\left|S_{\gamma}(f_{2})(x)\right|
&\lsim&\left[\int_{\Gamma(x)} \left(t^{-n}\int_{(2B)^{c}\cap\tilde{B}(u,t)}\left|f(z)\right|dz\right)^{2}\frac{dudt}{t^{n+1}}\right]^{\frac{1}{2}}\\
&\lsim&\sum^{\infty}_{k=1}\int_{2^{k+1}B\setminus 2^{k}B}\left|f(z)\right|\left[\int^{\infty}_{0}\left(\int_{B(x,t)}\chi_{\tilde{B}(z,t)}(u)du\right)\frac{dt}{t^{3n+1}}\right]^{\frac{1}{2}}dz
\end{eqnarray*}
where the last control is an application of Minkowski's integral inequality. 

We suppose $x\in B(y,r)$. For $k\in\mathbb N^{\ast}$, $z\in 2^{k+1}B\setminus 2^{k}B$ and $t>0$, $\int_{B(x,t)}\chi_{\tilde{B}(z,t)}(u)du\neq 0$ implies that $ B(x,t)\cap \tilde{B}(z,t)\neq\emptyset$. Let  $u_{0}\in B(x,t)\cap \tilde{B}(z,t)$, we have  
\begin{equation}
2t\geq \left|x-u_{0}\right|+\left|z-u_{0}\right|\geq\left|x-z\right|\geq \left|y-z\right|-\left|x-y\right|\geq 2^{k-1}r.\label{controlt}
\end{equation}
Thus for $x\in B=B(y,r)$,
\begin{eqnarray*}
\left|S_{\gamma}(f_{2})(x)\right|&\lsim&\sum^{\infty}_{k=1}\int_{2^{k+1}B\setminus 2^{k}B}\left|f(z)\right|\left(\int^{\infty}_{2^{k-2}r}\int_{B(x,t)}du\frac{dt}{t^{3n+1}}\right)^{\frac{1}{2}}dz\\
&\lsim&\sum^{\infty}_{k=1}\int_{2^{k+1}B\setminus 2^{k}B}\left|f(z)\right|\left(\int^{\infty}_{2^{k-2}r}\frac{dt}{t^{2n+1}}\right)^{\frac{1}{2}}dz\lsim\sum^{\infty}_{k=1}\frac{1}{\left|2^{k+1}B\right|}\int_{2^{k+1}B\setminus 2^{k}B}\left|f(z)\right|dz.
\end{eqnarray*}
But then by H\"older Inequality and (\ref{apcondition}), we have for every $k\in\mathbb N^{\ast}$
\begin{equation}
\frac{1}{\left|2^{k+1}B\right|}\int_{2^{k+1}B}\left|f(z)\right|dz\lsim
\left\|f\chi_{2^{k+1}B}\right\|_{q_{w}}w(2^{k+1}B)^{-\frac{1}{q}}.\label{contmoy}
\end{equation}
 It follows that 
\begin{equation}
\left\|S_{\gamma}(f_{2})\chi_{B(y,r)}\right\|_{q_{w}}\lsim\sum^{\infty}_{k=1} \left\|f\chi_{2^{k+1}B}\right\|_{q_{w}}\left(\frac{w(B)}{w(2^{k+1}B)}\right)^{\frac{1}{q}}.\label{estf2}
\end{equation}
Multiplying both Inequalities (\ref{estf1}) and (\ref{estf2}) by $w(B(y,r))^{\frac{1}{\alpha}-\frac{1}{q}-\frac{1}{p}}$, it comes from (\ref{contwl}) that 
\begin{equation}
\begin{aligned}
&w(B(y,r))^{\frac{1}{\alpha}-\frac{1}{q}-\frac{1}{p}}\left\|S_{\gamma}(f)\chi_{B(y,r)}\right\|_{q_{w}}\lsim w(B(y,2r))^{\frac{1}{\alpha}-\frac{1}{q}-\frac{1}{p}}\left\|f\chi_{B(y,2r)}\right\|_{q_{w}}\\
&\ \ \ \ \ \ \ \ \ \ \ \ \ \ \ +\sum^{\infty}_{k=1}w(B(y,2^{k+1}r))^{\frac{1}{\alpha}-\frac{1}{q}-\frac{1}{p}}\left\|f\chi_{B(y,2^{k+1}r)}\right\|_{q_{w}}\frac{1}{2^{\frac{nk}{s}(\frac{1}{\alpha}-\frac{1}{p})}},
\end{aligned}\label{3}
\end{equation}
 for some $s>0$. Therefore 
the $L^{p}$ norm of both sides of (\ref{3}) led to
\begin{equation*}
_{r}\left\|S_{\gamma}(f)\right\|_{q_{w},p,\alpha}\lsim(1+\sum^{\infty}_{k=1}\frac{1}{2^{\frac{nk}{s}(\frac{1}{\alpha}-\frac{1}{p})}}) \left\|f\right\|_{q_{w},p,\alpha},\ r>0,
\end{equation*}
and the result follows, since the series on the right hand side converge.
\epf

For the proof of Theorem \ref{main3}, we will need the following varying-aperture versions of $S_{\gamma}$. For $0<\gamma\leq 1$ and $\beta>0$, we define $S_{\gamma,\beta}(f)$ by
\begin{equation}
S_{\gamma,\beta}(f)(x)=\left[\int_{\Gamma_{\beta}(x)}\left(\sup_{\varphi\in\mathcal C_{\gamma}}\left|f\ast\varphi_{t}(y)\right|\right)^{2}\frac{dydt}{t^{n+1}}\right]^{\frac{1}{2}},
\end{equation}
where $\Gamma_{\beta}(x)=\left\{(x,t)\in\mathbb R^{n+1}_{+}/\left|x-y\right|<\beta t\right\}$. We have the following lemma which is a consequence of Lemmas 1.1, 1.2 and 1.3 of \cite{Wa} and the boundedness of $S_{\gamma}:=S_{\gamma,2^{0}}$ on the weighted Lebesgue spaces.
\begin{lem}\label{contsqw}Let $0<\gamma\leq 1$, $1<q<\infty$ and $w\in\mathcal A_{q}$. Then for all non negative integers $j$, $S_{\gamma,2^{j}}$ is bounded on $L^{q}_{w}(\mathbb R^{n})$. Moreover 
\begin{equation}
\left\|S_{\gamma,2^{j}}(f)\right\|_{q_{w}}\lsim (2^{nj}+2^{\frac{njq}{2}})\left\|f\right\|_{q_{w}}.
\end{equation} 
\end{lem}
 
\proof[Proof of Theorem \ref{main3}]
For all $x\in\mathbb R^{n}$, we have 
\begin{equation}
g^{\ast}_{\lambda,\gamma}(f)(x)^{2}\lsim S_{\gamma}(f)(x)^{2}+\sum^{\infty}_{j=1}2^{-j\lambda n}S_{\gamma,2^{j}}(f)(x)^{2}.\label{estg}
\end{equation}
Fix $r>0$. For $y\in \mathbb R^{n}$ and $B=B(y,r)$ a ball in $\mathbb R^{n}$, we have 
\begin{equation}
\begin{array}{lll}
w(B)^{\frac{1}{\alpha}-\frac{1}{q}-\frac{1}{p}}\left\|g^{\ast}_{\lambda,\gamma}(f)\chi_{B}\right\|_{q_{w}}&\lsim& w(B)^{\frac{1}{\alpha}-\frac{1}{q}-\frac{1}{p}}\left\|S_{\gamma}(f)\chi_{B}\right\|_{q_{w}}\\
&+&\sum^{\infty}_{j=1}2^{-j\lambda n/2}w(B)^{\frac{1}{\alpha}-\frac{1}{q}-\frac{1}{p}}\left\|S_{\gamma,2^{j}}(f)\chi_{B}\right\|_{q_{w}},
\end{array}\label{festimate}
\end{equation}
according to (\ref{estg}). By Theorem \ref{main1}, we have that the $L^{p}$ norm of the first term of (\ref{festimate}) is controlled by $\left\|f\right\|_{q_{w},p,\alpha}$. Let  $j$ be fixed in $\mathbb N^{\ast}$. For $S_{\gamma,2^{j}}f$, we proceed as for $S_{\gamma}$. So, for $f=f_{1}+f_{2}$ with $f_{1}=f\chi_{2B}$, we have 
\begin{equation}
\begin{array}{lll}w(B)^{\frac{1}{\alpha}-\frac{1}{q}-\frac{1}{p}}\left\|S_{\gamma,2^{j}}(f)\chi_{B}\right\|_{q_{w}}&\leq& w(B)^{\frac{1}{\alpha}-\frac{1}{q}-\frac{1}{p}}\left\|S_{\gamma,2^{j}}(f_{1})\chi_{B}\right\|_{q_{w}}\\
&+&w(B)^{\frac{1}{\alpha}-\frac{1}{q}-\frac{1}{p}}\left\|S_{\gamma,2^{j}}(f_{2})\chi_{B}\right\|_{q_{w}}.\label{estlp1}
\end{array}
\end{equation}
Applying Lemma \ref{contsqw} and taking into consideration  (\ref{doubling}), we obtain 
\begin{equation}
w(B)^{\frac{1}{\alpha}-\frac{1}{q}-\frac{1}{p}}\left\|S_{\gamma,2^{j}}(f_{1})\chi_{B}\right\|_{q_{w}}\lsim (2^{jn}+2^{jnq/2})w(2B)^{\frac{1}{\alpha}-\frac{1}{q}-\frac{1}{p}}\left\|f\chi_{2B}\right\|_{q_{w}}.\label{estlpf1}
\end{equation}
Let us estimate now the term in $f_{2}$. 
The same arguments we use to estimate $S_{\gamma}(f_{2})(x)$ for $x\in B$, i.e., Minkowsky's integral inequality and the fact that for $k\in\mathbb N^{\ast}$, $z\in 2^{k+1}B\setminus 2^{k}B$
$$\int_{B(x,2^{j}t)}\chi_{\tilde{B}(z,t)}(u)du\neq 0\Rightarrow t\geq\frac{2^{k-1}}{2^{j}+1}r,$$
 allow us to get the following 
\begin{equation*}
\begin{array}{lll}
\left|S_{\gamma,2^{j}}(f_{2})(x)\right|&\lsim& 2^{3jn/2}\sum^{\infty}_{k=1}\frac{1}{\left|2^{k+1}B\right|}\int_{2^{k+1}B\setminus 2^{k}B}\left|f(z)\right|dz\\
&\lsim& 2^{3jn/2}\sum^{\infty}_{k=1}\left\|f\chi_{B(y,2^{k+1}r)}\right\|_{q_{w}}w(B(y,2^{k+1}r))^{-\frac{1}{q}},
\end{array}
\end{equation*}
for all $x\in B(y,r)$, where the last control comes from  Estimation (\ref{contmoy}). Therefore, its $L^{q}_{w}(B)$-norm led to
\begin{equation}
w(B)^{\frac{1}{\alpha}-\frac{1}{q}-\frac{1}{p}}\left\|S_{\gamma,2^{j}}(f_{2})\chi_{B}\right\|_{q_{w}}\lsim2^{3jn/2}\sum^{\infty}_{k=1}\frac{w(2^{k+1}B)^{\frac{1}{\alpha}-\frac{1}{q}-\frac{1}{p}}}{2^{\frac{nk}{s}(\frac{1}{\alpha}-\frac{1}{p})}}\left\|f\chi_{2^{k+1}B}\right\|_{q_{w}}.\label{estlpf2}
\end{equation}
Taking Estimates (\ref{estlpf1}) and (\ref{estlpf2}) in (\ref{estlp1}), we have
\begin{equation*}
\begin{array}{lll}
w(B(y,r))^{\frac{1}{\alpha}-\frac{1}{q}-\frac{1}{p}}\left\|S_{\gamma,2^{j}}(f)\chi_{B(y,r)}\right\|_{q_{w}}&\lsim&(2^{jn}+2^{jnq/2})w(B(y,2r))^{\frac{1}{\alpha}-\frac{1}{q}-\frac{1}{p}}\left\|f\chi_{B(y,2r)}\right\|_{q_{w}}\\
&+&2^{3jn/2}\sum^{\infty}_{k=1}\frac{w(B(y,2^{k+1}r))^{\frac{1}{\alpha}-\frac{1}{q}-\frac{1}{p}}}{2^{\frac{nk}{s}(\frac{1}{\alpha}-\frac{1}{p})}}\left\|f\chi_{B(y,2^{k+1}r)}\right\|_{q_{w}}
\end{array},
\end{equation*}
for all $y\in\mathbb R^{n}$, so that the $L^{p}$-norm of both sides led to
\begin{equation}
\ _{r}\left\|S_{\gamma,2^{j}}(f)\right\|_{q_{w},p,\alpha}\lsim (2^{jn}+2^{jnq/2})\left\|f\right\|_{q_{w},p,\alpha}+  \left\|f\right\|_{q_{w},p,\alpha}2^{3jn/2}. 
\end{equation}
Therefore the $L^{p}$ norm of (\ref{festimate}) give
\begin{equation}
\begin{array}{lll}
\ _{r}\left\|g^{\ast}_{\lambda,\gamma}(f)\right\|_{q_{w},p,\alpha}&\lsim&\left(1+\sum^{\infty}_{j=1}2^{-j\lambda n/2}(2^{jn}+2^{jnq/2}+2^{3jn/2})\right)\left\|f\right\|_{q_{w},p,\alpha}\\
&\lsim& \left\|f\right\|_{q_{w},p,\alpha}, 
\end{array}
\end{equation}
for $r>0$, where the convergence of the series is due to the fact that $\lambda>\max\left\{q,3\right\}$. We end the proof by taking the supremum over all $r>0$.
\epf

For the proof of the next results on commutators, we use the following properties of $BMO$ (see \cite{JN}). Let $b$ be a locally integrable function. If $b\in BMO(\mathbb R^{n})$, then for every $1<p<\infty$, we have
\begin{equation}
\left\|b\right\|_{BMO(\mathbb R^{n})}\approx\sup_{B:\text{ ball}}\left(\frac{1}{\left|B\right|}\int_{B}\left|b(x)-b_{B}\right|^{p}dx\right)^{\frac{1}{p}},\label{equivbmo}
\end{equation}
and for $w\in\mathcal A_{q}$ with $1<q<\infty$,
\begin{equation}
\left(\frac{1}{w(B)}\int_{B}\left|b(x)-b_{B}\right|^{p}w(x)dx\right)^{\frac{1}{p}}\lsim \left\|b\right\|_{BMO},\label{controlw}
\end{equation}
which is an immediate consequence of (\ref{equivbmo}) and the characterization (\ref{reverseholder}) of $\mathcal A_{q}$ weights.
\proof[Proof of Theorem \ref{main2}]
Fix $y\in\mathbb R^{n}$ and $r>0$. For $B=B(y,r)$, we put $f=f_{1}+f_{2}$ with $f_{1}=f\chi_{2B}$. We have 
\begin{equation}
\begin{array}{lll}
w(B)^{\frac{1}{\alpha}-\frac{1}{q}-\frac{1}{p}}\left\|\left[b,S_{\gamma}\right](f)\chi_{B}\right\|_{q_{w}}&\leq& w(B)^{\frac{1}{\alpha}-\frac{1}{q}-\frac{1}{p}}\left\|\left[b,S_{\gamma}\right](f_{1})\chi_{B}\right\|_{q_{w}}\\
&+&w(B)^{\frac{1}{\alpha}-\frac{1}{q}-\frac{1}{p}}\left\|\left[b,S_{\gamma}\right](f_{2})\chi_{B}\right\|_{q_{w}}
\end{array}\label{estcomm1}
\end{equation}
For the term in $f_{1}$, it is immediate that  
\begin{equation}
w(B(y,r))^{\frac{1}{\alpha}-\frac{1}{q}-\frac{1}{p}}\left\|\left[b,S_{\gamma}\right](f_{1})\chi_{B(y,r)}\right\|_{q_{w}}
\lsim w(B(y,2r))^{\frac{1}{\alpha}-\frac{1}{q}-\frac{1}{p}}\left\|f\chi_{B(y,2r)}\right\|_{q_{w}},\label{estcommf1}
\end{equation}
 according to the boundedness of the commutator on $L^{q}_{w}(\mathbb R^{n})$ and (\ref{doubling}). It remains to estimate the term in $f_{2}$.

Let $x\in\mathbb R^{n}$. For $(u,t)\in\Gamma(x)$, we have
\begin{equation*}
\begin{aligned}
&\sup_{\varphi\in\mathcal C_{\gamma}}\left|\int_{\mathbb R^{n}}(b(x)-b(z))\varphi_{t}(u-z)f_{2}(z)dz\right|\\
&\ \ \ \ \ \ \ \ \ \ \ \ \ \ \ \ \ \ \ \leq\left|b(x)-b_{B}\right| \sup_{\varphi\in\mathcal C_{\gamma}}\left|f_{2}\ast\varphi_{t}(u)\right|
+\sup_{\varphi\in\mathcal C_{\gamma}}\left|(b-b_{B})f_{2}\ast\varphi_{t}(u)\right| 
\end{aligned}
\end{equation*}
so that the $L^{2}(\Gamma(x),\frac{dudt}{t^{n+1}})$-norm of both sides led to
\begin{equation*}
\begin{array}{lll}
\left|\left[b,S_{\gamma}\right]f_{2}(x)\right|&\leq& \left|b(x)-b_{B}\right|S_{\gamma}(f_{2})(x)\\
&+&\left\{\int_{\Gamma(x)}\left(\sup_{\varphi\in\mathcal C_{\gamma}}\left|[(b-b_{B})f_{2}]\ast\varphi_{t}(u)\right|\right)^{2}
\frac{dudt}{t^{n+1}}\right\}^{\frac{1}{2}}=I+II.
\end{array}
\end{equation*}
We take $x\in B=B(y,r)$. As we prove in Theorem \ref{main1}, we have  
$$\left|S_{\gamma}(f_{2})(x)\right|\lsim\sum^{\infty}_{k=1} \left\|f\chi_{2^{k+1}B}\right\|_{q_{w}}w(2^{k+1}B)^{-\frac{1}{q}}.$$
 Thus, the $L^{q}_{w}(B)$-norm of $I$ can be estimated as follow 
\begin{equation}
\left\|\left|b-b_{B}\right|S_{\gamma}(f_{2})\chi_{B}\right\|_{q_{w}} \lsim\left\|b\right\|_{BMO}\sum^{\infty}_{k=1}\left(\frac{w(B)}{w(2^{k+1}B)}\right)^{\frac{1}{q}}\left\|f\chi_{2^{k+1}B}\right\|_{q_{w}},\label{estimationI}
\end{equation}
where we use (\ref{controlw}). 
 On other hand, it comes from the uniformly boundedness of the family $\mathcal C_{\gamma}$ that 
\begin{equation*}
II\lsim \left[\int_{\Gamma(x)}\left(t^{-n}\int_{(2B)^{c}\cap\tilde{B}(u,t)}\left|b(z)-b_{B}\right|\left|f(z)\right|dz\right)^{2}\frac{dudt}{t^{n+1}}\right]^{\frac{1}{2}},
\end{equation*}
so that using once more Minkowski's inequality for integrals and Inequality (\ref{controlt}), we have 
\begin{eqnarray*}
II&\lsim&\sum^{\infty}_{k=1}\frac{1}{\left|2^{k+1}B\right|}\int_{2^{k+1}B\setminus 2^{k}B}\left|b(z)-b_{B}\right|\left|f(z)\right|dz\\
&\leq&\sum^{\infty}_{k=1}\frac{1}{\left|2^{k+1}B\right|}\int_{2^{k+1}B\setminus 2^{k}B}\left|b(z)-b_{2^{k+1}B}\right|\left|f(z)\right|dz+\sum^{\infty}_{k=1}\frac{\left|b_{2^{k+1}B}-b_{B}\right|}{\left|2^{k+1}B\right|}\int_{2^{k+1}B\setminus 2^{k}B}\left|f(z)\right|dz
\end{eqnarray*} 
for all $x\in B(y,r)$. 
But then 
\begin{eqnarray*}
\int_{(2^{k+1}B\setminus 2^{k}B)}\left|b(z)-b_{2^{k+1}B}\right|\left|f(z)\right|dz&\leq&\left(\int_{2^{k+1}B}\left|b(z)-b_{2^{k+1}B}\right|^{q'}w(z)^{-\frac{q'}{q}}dz\right)^{\frac{1}{q'}}\\
&\times&\left(\int_{2^{k+1}B}\left|f(z)\right|^{q}w(z)dz\right)^{\frac{1}{q}}\\
&\lsim&\left\|f\chi_{2^{k+1}B}\right\|_{q_{w}}\left|2^{k+1}B\right|w(2^{k+1}B)^{-\frac{1}{q}}\left\|b\right\|_{BMO},\label{estbmo}
%
\end{eqnarray*}
according to H\"older inequality and the fact that the weight $v(z)=w(z)^{-\frac{q'}{q}}$ belongs to $\mathcal A_{q'}$ whenever $w\in\mathcal A_{q}$. So
\begin{eqnarray*}
\sigma_{1}&:=&\sum^{\infty}_{k=1}\frac{1}{\left|2^{k+1}B\right|}\int_{2^{k+1}B\setminus 2^{k}B}\left|b(z)-b_{2^{k+1}B}\right|\left|f(z)\right|dz\\
&\lsim& \left\|b\right\|_{BMO}\sum^{\infty}_{k=1}\left\|f\chi_{2^{k+1}B}\right\|_{q_{w}}w(2^{k+1}B)^{-\frac{1}{q}},
\end{eqnarray*}
on $B$, and 
 the $L^{q}_{w}(B)$ norm of both sides led to 
\begin{equation}
\left\|\sigma_{1}\chi_{B(y,r)}\right\|_{q_{w}}\lsim\left\|b\right\|_{BMO}\sum^{\infty}_{k=1}\left\|f\chi_{2^{k+1}B}\right\|_{q_{w}}\left(\frac{w(B)}{w(2^{k+1}B)}\right)^{\frac{1}{q}}.\label{sigma1}
\end{equation}
For the second series, we have
\begin{eqnarray*}
\sigma_{2}&:=&\sum^{\infty}_{k=1}\frac{\left|b_{2^{k+1}B}-b_{B}\right|}{\left|2^{k+1}B\right|}\int_{2^{k+1}B\setminus 2^{k}B}\left|f(z)\right|dz\\
&\lsim& \left\|b\right\|_{BMO}\left(\sum^{\infty}_{k=1}(k+1)\left\|f\chi_{2^{k+1}B}\right\|_{q_{w}}w(2^{k+1}B)^{-\frac{1}{q}}\right),
\end{eqnarray*}
where we use the fact that 
$\left|b_{2^{k+1}B}-b_{B}\right|\lsim (k+1)\left\|b\right\|_{BMO}$ and Relation (\ref{contmoy}).
It comes that 
\begin{equation}
\left\|\sigma_{2}\chi_{B(y,r)}\right\|_{q_{w}}\lsim\left\|b\right\|_{BMO}\left(\sum^{\infty}_{k=1}(k+1)\left\|f\chi_{2^{k+1}B}\right\|_{q_{w}}\left(\frac{w(B)}{w(2^{k+1}B)}\right)^{\frac{1}{q}}\right)\label{sigma2}.
\end{equation}
Hence, putting together (\ref{estimationI}), (\ref{sigma1}) and (\ref{sigma2}), we obtain, 
\begin{equation}
\begin{aligned}
&w(B(y,r))^{\frac{1}{\alpha}-\frac{1}{q}-\frac{1}{p}}\left\|\left[b,S_{\eta}\right](f_{2})\chi_{B(y,r)}\right\|_{q_{w}}\\
&\ \ \lsim \left\|b\right\|_{BMO}\left(\sum^{\infty}_{k=1}\frac{k+3}{2^{\frac{2nk}{s'}(\frac{1}{\alpha}-\frac{1}{p})}}w(B(y,2^{k+1}r))^{\frac{1}{\alpha}-\frac{1}{q}-\frac{1}{p}}\left\|f\chi_{B(y,2^{k+1}r)}\right\|_{q_{w}}\right)
\end{aligned}\label{estcommf2}
\end{equation}
for all $y\in\mathbb R^{n}$ and some $s'> 0$.
Taking Estimates (\ref{estcommf1}) and (\ref{estcommf2}) in (\ref{estcomm1}) yield,
\begin{equation}
\begin{aligned}
&w(B(y,r))^{\frac{1}{\alpha}-\frac{1}{q}-\frac{1}{p}}\left\|\left[b,S_{\eta}\right](f)\chi_{B(y,r)}\right\|_{q_{w}}\\
&\ \ \ \ \ \ \lsim\left\|b\right\|_{BMO}\left(\sum^{\infty}_{k=1}\frac{k+3}{2^{\frac{2nk}{s'}(\frac{1}{\alpha}-\frac{1}{p})}}w(B(y,2^{k+1}r))^{\frac{1}{\alpha}-\frac{1}{q}-\frac{1}{p}}\left\|f\chi_{B(y,2^{k+1}r)}\right\|_{q_{w}}\right)\\
&\ \ \ \ \ \ + w(B(y,2r))^{\frac{1}{\alpha}-\frac{1}{q}-\frac{1}{p}}\left\|f\chi_{B(y,2r)}\right\|_{q_{w}}
\end{aligned}\label{estcomm2}
\end{equation}
for all $y\in\mathbb R^{n}$. Therefore the $L^{p}$-norm of both sides of (\ref{estcomm2}), gives 
$$\ _{r}\left\|\left[b,S_{\gamma}\right](f)\right\|_{q_{w},p,\alpha}\lsim(1+\left\|b\right\|_{BMO})\left\|f\right\|_{q_{w},p,\alpha},$$
for all $r>0$, since the series $\sum^{\infty}_{k=1}\frac{k+3}{2^{\frac{2nk}{s'}(\frac{1}{\alpha}-\frac{1}{p})}}$ converges. We end the proof by taking the supremum over all $r>0$.
\epf

\proof[Proof of Theorem \ref{nouveau}]
It is easy to see that 
\begin{equation*}
\left[b,g^{\ast}_{\lambda,\gamma}\right](f)^{2}(x)\lsim \sum^{\infty}_{j=0}2^{-j\lambda n}\left[b,S_{\gamma,2^{j}}\right](f)^{2}(x),
\end{equation*}
for all $x\in\mathbb R^{n}$. So, for all balls $B=B(y,r)$ we have
\begin{equation*}
\left\|\left[b,g^{\ast}_{\lambda,\gamma}\right](f)\chi_{B}\right\|_{q_{w}}\lsim \sum^{\infty}_{j=0}2^{-\frac{j\lambda n}{2}}\left\|\left[b,S_{\gamma,2^{j}}\right](f)\chi_{B}\right\|_{q_{w}}.
\end{equation*}
Using the arguments as in the proof of theorems \ref{main2} and \ref{main1} and taking into consideration (\ref{estlpf1}) we end the proof.
\epf

\end{document}